\documentclass[graybox]{svmult}

\usepackage{amsmath}
\usepackage{amssymb}
\usepackage{tabu}
\usepackage[caption=false]{subfig}
\usepackage{tikz}
\usepackage{url}

\usepackage{type1cm}        
\usepackage{makeidx}         
\usepackage{graphicx}        
\usepackage{multicol}        
\usepackage[bottom]{footmisc}

\usepackage{newtxtext}       %
\usepackage{newtxmath}       

\makeindex             

\begin{document}

\title*{On the Dirichlet-to-Neumann coarse space for solving the Helmholtz problem using domain decomposition}
\titlerunning{On the Dirichlet-to-Neumann coarse space for the Helmholtz problem}
\author{Niall Bootland and Victorita Dolean}
\institute{Niall Bootland \at University of Strathclyde, Department of Mathematics and Statistics, Glasgow, UK.\\ \email{niall.bootland@strath.ac.uk}
\and Victorita Dolean \at University of Strathclyde, Deptartment of Mathematics and Statistics, Glasgow, UK.\\ Universit\'e C\^ote d'Azur, CNRS, Laboratoire J.A. Dieudonn\'e, Nice, France.\\ \email{work@victoritadolean.com}}
%
%
\maketitle

\abstract*{We examine the use of the Dirichlet-to-Neumann coarse space within an additive Schwarz method to solve the Helmholtz equation in 2D. In particular, we focus on the selection of how many eigenfunctions should go into the coarse space. We find that wave number independent convergence of a preconditioned iterative method can be achieved in certain special cases with an appropriate and novel choice of threshold in the selection criteria. However, this property is lost in a more general setting, including the heterogeneous problem. Nonetheless, the approach converges in a small number of iterations for the homogeneous problem even for relatively large wave numbers and is robust to the number of subdomains used.}

\abstract{We examine the use of the Dirichlet-to-Neumann coarse space within an additive Schwarz method to solve the Helmholtz equation in 2D. In particular, we focus on the selection of how many eigenfunctions should go into the coarse space. We find that wave number independent convergence of a preconditioned iterative method can be achieved in certain special cases with an appropriate and novel choice of threshold in the selection criteria. However, this property is lost in a more general setting, including the heterogeneous problem. Nonetheless, the approach converges in a small number of iterations for the homogeneous problem even for relatively large wave numbers and is robust to the number of subdomains used.}

\section{Introduction}
\label{sec:1-introduction}
Within domain decomposition methods, the use of a coarse space as a second level is typically required to provide scalability with respect to the number of subdomains used \cite{DoleanEtAl-Book}. More recently, coarse spaces have also been designed to provide robustness to model parameters, especially for large contrasts in heterogeneous problems. For example, the GenEO coarse space has been successfully employed for the robust solution of highly heterogeneous elliptic problems \cite{SpillaneEtAl-GenEO}. One way in which a coarse space can be derived is via solving local eigenvalue problems on subdomains, as is the case for the GenEO method. An earlier approach, having many similarities, is the Dirichlet-to-Neumann (DtN) coarse space \cite{NatafEtAl-DtN}. We focus on this method which solves eigenvalue problems on the boundary of subdomains related to a Dirichlet-to-Neumann map.

We are interested in using domain decomposition methodology to solve wave propagation problems. In particular, we consider the Helmholtz problem\footnote{Note that if $\Gamma_{R} = \emptyset$ then the problem will be ill-posed for certain choices of $k$ corresponding to Dirichlet eigenvalues of the corresponding Laplace problem.}\begin{subequations}
\label{Helmholtz}
\begin{align}
\label{HelmholtzEquation}
-\Delta u - k^{2} u & = f & & \text{in } \Omega,\\
\label{HelmholtzDirichletBC}
u & = 0 & & \text{on } \Gamma_{D},\\
\label{HelmholtzRobinBC}
\frac{\partial u}{\partial n} + i k u & = 0 & & \text{on } \Gamma_{R},
\end{align}
\end{subequations}
with wave number $k > 0$, where $\partial\Omega = \Gamma_{D}\cup\Gamma_{R}$ and $\Gamma_{D}\cap\Gamma_{R} = \emptyset$. Such problems arise in many wave propagation and scattering problems in science and engineering, for instance, acoustic and seismic imaging problems. Furthermore, we also consider the heterogeneous problem, in which case $k(\vec{x})$ varies in the domain $\Omega$. We suppose the variation in $k$ stems from the wave speed $c(\vec{x})$ depending on the heterogeneous media, with the wave number being given by $k = \omega/c$ for angular frequency $\omega$.

The wave number $k$ is the key parameter within the Helmholtz equation and as $k$ increases the problem becomes more challenging. We are interested in the case when $k$ becomes large and so solutions are highly oscillatory. The numerical method employed needs to be able to capture this behaviour, often through an increasing number of grid points, such as a fixed number of points per wavelength. However, typically the number of grid points needs to grow faster than linearly in $k$ if accuracy is to be maintained due to the pollution effect \cite{BabuskaAndSauter-Pollution}. For instance, when using P1 finite elements for the numerical solution of \eqref{Helmholtz}, the mesh spacing $h$ should grow as $k^{-3/2}$. This means very large linear systems must be solved when $k$ is large and, since these systems are sparse, iterative methods are most often employed for their solution. However, efficiently solving large discrete Helmholtz systems is challenging since classical iterative methods fail to be effective \cite{ErnstAndGander-Helmholtz}. As such, we require a more robust iterative solver. Here we consider a restricted additive Schwarz (RAS) method with a Dirichlet-to-Neumann coarse space \cite{ConenEtAl-DtN} and will be interested in the performance of this solver methodology as $k$ increases. We now review the underlying numerical methods we use.

\section{Discretisation and solver methodology}
\label{sec:2-DiscretisationAndSolverMethodology}

To discretise we use finite element methodology, in particular using piecewise linear (P1) finite elements on simplicial meshes. Given a simplicial mesh $T^{h}$ on a bounded polygonal domain $\Omega$, let $V^{h} \subset \left\lbrace H^{1}(\Omega) \colon u = 0 \text{ on } \Gamma_{D} \right\rbrace$ be the space of piecewise linear functions on $T^{h}$. The P1 finite element solution $u_{h} \in V^{h}$ satisfies the weak formulation $a(u_{h},v_{h}) = F(v_{h}) \ \forall v_{h} \in V^{h}$, where
\begin{align}
\label{WeakForms}
a(u,v) & = \int_{\Omega} \left( \nabla u \cdot \nabla \bar{v} - k^2 u \bar{v}\right) \, \mathrm{d}\vec{x} + \int_{\Gamma_{R}} i k u \bar{v} \, \mathrm{d}s, & \text{and} & & F(v) & = \int_{\Omega} f \bar{v} \, \mathrm{d}\vec{x}.
\end{align}
Using the standard nodal basis for $V^{h}$ we can represent the solution $u_{h}$ through its basis coefficients $\vec{u}$ and reduce the problem to solving the complex symmetric linear system $A\vec{u} = \vec{f}$ where $A$ comes from the bilinear form $a(\cdot,\cdot)$ and $\vec{f}$ the linear functional $F(\cdot)$; see, for example, \cite{ConenEtAl-DtN}.

To solve the discrete Helmholtz system $A\vec{u} = \vec{f}$ we utilise a two-level domain decomposition preconditioner within an iterative Krylov method. Since $A$ is only complex symmetric rather than Hermitian, we use GMRES as the iterative Krylov method \cite{SaadAndSchultz-GMRES}. For the domain decomposition, given an overlapping partition $\left\lbrace \Omega_{j} \right\rbrace_{j=1}^{N}$ of $\Omega$, let $R_{j}$ represent the matrix form of the restriction onto subdomain $\Omega_{j}$. Then the restricted additive Schwarz (RAS) domain decomposition preconditioner is given by
\begin{align}
\label{RAS}
M_{\text{RAS}}^{-1} = \sum_{j=1}^{N} R_{j}^{T} D_{j} A_{j}^{-1} R_{j},
\end{align}
where $A_{j} = R_{j}AR_{j}^{T}$ is the local Dirichlet matrix on $\Omega_{j}$ and the diagonal matrices $D_{j}$ are a discrete representation of a partition of unity (see \cite{DoleanEtAl-Book}); this removes ``double counting''  in regions of overlap. Note that each subdomain contribution from the sum in \eqref{RAS} can be computed locally in parallel. Using the one-level preconditioner \eqref{RAS} is not sufficient to provide robustness with respect to the number of subdomains $N$ used and also becomes much worse when $k$ increases. To this end we incorporate a coarse space as a second level within the method.

A coarse space provides a more efficient way to transfer information globally between subdomains, rather than relying solely on local solutions, as in \eqref{RAS}. The coarse space constitutes a collection of column vectors $Z$, having full column rank. We then utilise the coarse correction operator $Q = Z E^{-1} Z^{\dagger}$, where $E = Z^{\dagger} A Z$ is the coarse space operator, which provides a coarse solution in the space spanned by the columns of $Z$. To incorporate the coarse correction we use an adapted deflation (AD) approach given by the two-level preconditioner
\begin{align}
\label{2LevelAdaptiveDeflationPreconditioner}
M_{AD}^{-1} = M_{\text{RAS}}^{-1}(I-AQ) + Q.
\end{align}
To complete the specification, we must choose which vectors go into the coarse space matrix $Z$.

\section{The Dirichlet-to-Neumann coarse space}
\label{sec:3-TheDirichetletToNeummanCoarseSpace}

We now introduce the Dirichlet-to-Neumann coarse space. The construction is based on solving local eigenvalue problems on subdomain boundaries related to the DtN map. To define this map we first require the Helmholtz extension operator from the boundary of a subdomain $\Omega_{j}$.

Let $\Gamma_{j} = \partial\Omega_{j} \setminus \partial\Omega$ and suppose we have Dirichlet data $v_{\Gamma_{j}}$ on $\Gamma_{j}$, then the Helmholtz extension $v$ in $\Omega_{j}$ is defined as the solution of\begin{subequations}
\label{HelmholtzExtension}
\begin{align}
\label{HelmholtzExtensionEquation}
-\Delta v - k^{2} v & = 0 & & \text{in } \Omega_{j},\\
\label{HelmholtzExtensionDirichletBC}
v & = v_{\Gamma_{j}} & & \text{on } \Gamma_{j},\\
\label{HelmholtzExtensionProblemBC}
\mathcal{C}(v) & = 0 & & \text{on } \partial\Omega_{j}\cap\partial\Omega,
\end{align}
\end{subequations}
where $\mathcal{C}(v) = 0$ represents the original problem boundary conditions \eqref{HelmholtzDirichletBC}--\eqref{HelmholtzRobinBC}. The DtN map takes Dirichlet data $v_{\Gamma_{j}}$ on $\Gamma_{j}$ to the corresponding Neumann data, that is
\begin{align}
\label{DtNMap}
\mathrm{DtN}_{\Omega_{j}}(v_{\Gamma_{j}}) = \left.\frac{\partial v}{\partial n} \right\rvert_{\Gamma_{j}}
\end{align}
where $v$ is the Helmholtz extension defined by \eqref{HelmholtzExtension}.

We now seek eigenfunctions of the DtN map locally on each subdomain $\Omega_{j}$, given by solving
\begin{align}
\label{DtNEigenproblem}
\mathrm{DtN}_{\Omega_{j}}(u_{\Gamma_{j}}) = \lambda u_{\Gamma_{j}}
\end{align}
for eigenfunctions $u_{\Gamma_{j}}$ and eigenvalues $\lambda \in \mathbb{C}$. To provide functions to go into the coarse space, we take the Helmholtz extension of $u_{\Gamma_{j}}$ in $\Omega_{j}$ and then extend by zero into the whole domain $\Omega$ using the partition of unity. For further details and motivation, as well as the discrete formulation of the eigenproblems, see \cite{ConenEtAl-DtN}.

It remains to determine which eigenfunctions of \eqref{DtNEigenproblem} should be included in the coarse space. Several selection criteria were investigated in \cite{ConenEtAl-DtN} and it was clear that the best choice was to select eigenvectors corresponding to eigenvalues with the smallest real part. That is, we use a threshold on the abscissa $\eta = \mathrm{Re}(\lambda)$ given by
\begin{align}
\label{Threshold}
\eta < \eta_{\text{max}}
\end{align}
where $\eta_{\text{max}}$ depends on $k_{j} = \max_{\vec{x}\in\Omega_{j}} k(\vec{x})$. In particular, \cite{ConenEtAl-DtN} advocate the choice $\eta_{\text{max}} = k_{j}$. Clearly, the larger $\eta_{\text{max}}$ is taken, the more eigenfunctions we include in the coarse space, increasing its size and the associated computational cost. However, it is not clear that $\eta_{\text{max}} = k_{j}$ is necessarily the best choice. We investigate the utility of choosing $\eta_{\text{max}}$ larger than $k_{j}$ and will see that, in some cases, taking a slightly larger coarse space can give improved behaviour as the wave number $k$ is increased.

\section{Numerical results}
\label{sec:4-NumericalResults}

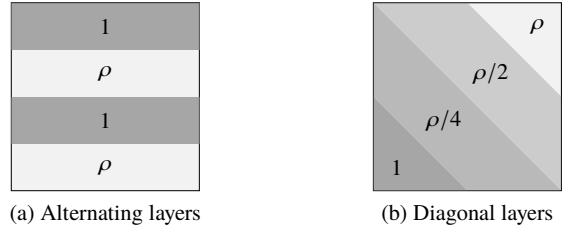
\begin{figure}[t]
	\centering
	\hspace*{\fill}
	\subfloat[Alternating layers]{
		\begin{tikzpicture}[scale=1.25]
		\draw (-1,1) -- (-1,-1) -- (1,-1) -- (1,1) -- cycle;
		\fill[gray,opacity=0.7] (-1,1) -- (1,1) -- (1,1/2) -- (-1,1/2) -- cycle;
		\fill[gray,opacity=0.1] (-1,0) -- (1,0) -- (1,1/2) -- (-1,1/2) -- cycle;
		\fill[gray,opacity=0.7] (-1,0) -- (1,0) -- (1,-1/2) -- (-1,-1/2) -- cycle;
		\fill[gray,opacity=0.1] (-1,-1) -- (1,-1) -- (1,-1/2) -- (-1,-1/2) -- cycle;
		\draw (0,3/4) node {$1$};
		\draw (0,1/4) node {$\rho$};
		\draw (0,-1/4) node {$1$};
		\draw (0,-3/4) node {$\rho$};
		\label{Fig:WaveGuideLayersA}
		\end{tikzpicture}
	} \hspace*{\fill}
	\subfloat[Diagonal layers]{
		\begin{tikzpicture}[scale=1.25]
		\draw (-1,1) -- (-1,-1) -- (1,-1) -- (1,1) -- cycle;
		\fill[gray,opacity=0.1] (-0,1) -- (1,1) -- (1,0) -- cycle;
		\fill[gray,opacity=0.4] (-1,1) -- (0,1) -- (1,0) -- (1,-1) -- cycle;
		\fill[gray,opacity=0.55] (-1,1) -- (-1,0) -- (0,-1) -- (1,-1) -- cycle;
		\fill[gray,opacity=0.7] (-1,0) -- (-1,-1) -- (0,-1) -- cycle;
		\draw (-3/4,-3/4) node {$1$};
		\draw (-1/4,-1/4) node {$\rho/4$};
		\draw (1/4,1/4) node {$\rho/2$};
		\draw (3/4,3/4) node {$\rho$};
		\label{Fig:WaveGuideLayersB}
		\end{tikzpicture}
	}
	\hspace*{\fill}
	\caption{Different layered configurations for the heterogeneous wave speed $c(\vec{x})$ within the wave guide problem, where $\rho>1$ is a contrast parameter.}
	\label{Fig:WaveGuideLayers}
\end{figure}

To investigate the dependence on $\eta_{\text{max}}$ we use a 2D wave guide problem on the unit square $\Omega = (0,1)^2$ as our model test problem. The Dirichlet condition \eqref{HelmholtzDirichletBC} is imposed on the left and right boundaries $\Gamma_{D} = \lbrace0,1\rbrace\times[0,1]$ while the Robin condition \eqref{HelmholtzRobinBC} is prescribed for the top and bottom boundaries $\Gamma_{R} = [0,1]\times\lbrace0,1\rbrace$. The right-hand side $f$ models a point source at the centre $(\frac{1}{2},\frac{1}{2})$. The wave number $k$ is either constant throughout $\Omega$ for the homogeneous problem or else $k = \omega/c$ where $\omega$ is constant and $c(\vec{x})$ is piecewise constant as illustrated in Figure \ref{Fig:WaveGuideLayers} for a contrast parameter $\rho > 1$. These heterogeneous problems model layered media.

To discretise we use a uniform square grid with $n_{\text{glob}}$ points in each direction and triangulate with alternating diagonals to form the P1 elements. As we increase $k$ we choose $n_{\text{glob}} \propto k^{3/2}$ in order to ameliorate the pollution effect. To begin with, we use a uniform decomposition into $N$ square subdomains and throughout use minimal overlap. All computations are performed using FreeFem (\url{http://freefem.org/}), in particular using the \texttt{ffddm} framework. When solving the linear systems we use preconditioned GMRES with the two-level preconditioner \eqref{2LevelAdaptiveDeflationPreconditioner} incorporating the DtN coarse space with threshold $\eta_{\text{max}}$ to reach a relative residual tolerance of $10^{-6}$.

\begin{figure}[t]
	\centering
	\hspace*{\fill}
	\subfloat{
		\includegraphics[width=0.43\linewidth,clip,trim=95mm 3mm 100mm 9mm]{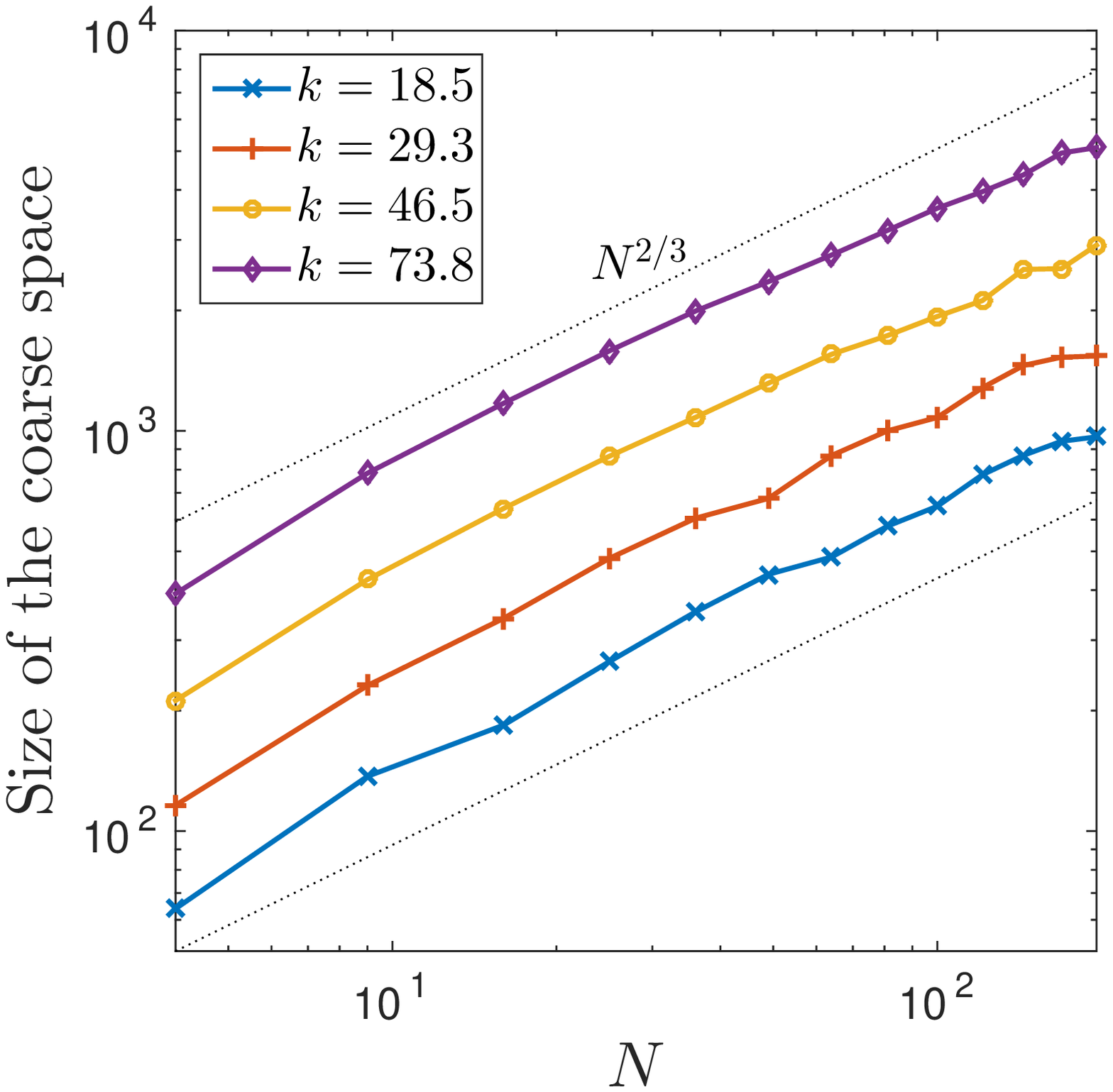}
	} \hspace*{\fill}
	\subfloat{
		\includegraphics[width=0.43\linewidth,clip,trim=95mm 3mm 100mm 9mm]{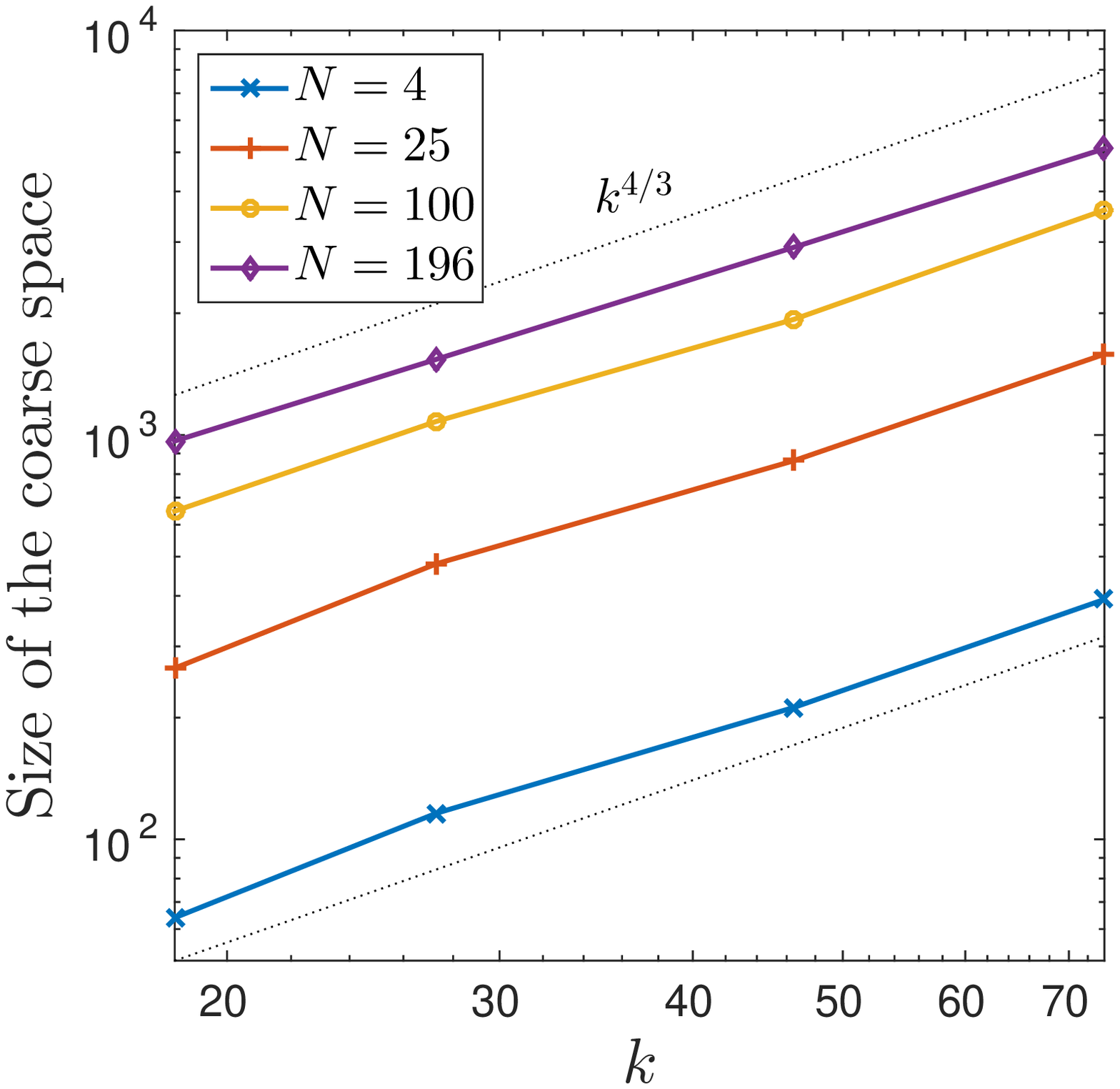}
	} \hspace*{\fill}
	\caption{The size of the DtN coarse space as a function of the number of subdomains $N$ (left) and wave number $k$ (right) for the homogeneous problem with threshold $\eta_{\text{max}}=k^{4/3}$.}
	\label{Fig:DependencePlots}
\end{figure}

In Table \ref{Table:HomoHelm5x5DtNVaryingThreshold} we vary the threshold $\eta_{\text{max}}$ as powers $k$ for the homogeneous problem using a fixed $5\times5$ square decomposition. The best choice advocated in \cite{ConenEtAl-DtN}, namely $\eta_{\text{max}} = k$, succeeds in requiring relatively low iteration counts to reach convergence with a modest size of coarse space. However, we observe that as the wave number $k$ increases the number of iterations required also increases, suggesting the approach will begin to struggle if $k$ becomes too large. We see from other choices of $\eta_{\text{max}}$ that taking a larger coarse space reduces the iteration counts. For instance, with the largest wave number tested when $\eta_{\text{max}} = k^{1.2}$ the size of the coarse space doubles while the iteration count it cut almost by a factor of three compared to $\eta_{\text{max}}=k$. In fact, there is a qualitative change in behaviour with respect to the wave number $k$, namely independence of the iteration counts to $k$, once $\eta_{\text{max}}$ becomes large enough, this point being approximately given by $\eta_{\text{max}}=k^{4/3}$. We note that the size of the coarse space is approximately proportional to $\eta_{\text{max}}$ in the results of Table \ref{Table:HomoHelm5x5DtNVaryingThreshold} (see also Figure \ref{Fig:DependencePlots}). As such, we see that the coarse space should grow faster than linearly in $k$ in order to achieve wave number independent iteration counts for this problem.

\begin{table}[ht]
	\centering
	\caption{Preconditioned GMRES iteration counts using the two-level method while varying the threshold parameter $\eta_{\text{max}}$ for the DtN coarse space. The size of the coarse space is given in brackets. A uniform decomposition into $5\times5$ square subdomains is used.}
	\label{Table:HomoHelm5x5DtNVaryingThreshold}
	\tabulinesep=1mm
	\begin{tabu}{cc|cccccc}
		$n_{\text{glob}}$ & $k$ & $\ \ \eta_{\text{max}} = k\ \ $ & $\ \eta_{\text{max}} = k^{1.1}\ $ & $\ \eta_{\text{max}} = k^{1.2}\ $ & $\ \eta_{\text{max}} = k^{1.3}\ $ & $\ \eta_{\text{max}} = k^{1.4}\ $ & $\ \eta_{\text{max}} = k^{1.5}\ $ \\
		\hline
		100 & 18.5 & 12 (144) &  9 (160) &  8 (200) & 7 (240)  & 6 (320)  & 5 (400)  \\
		200 & 29.3 & 16 (215) & 11 (240) &  9 (320) & 7 (434)  & 6 (560)  & 5 (760)  \\
		400 & 46.5 & 18 (299) & 13 (393) & 10 (545) & 7 (784)  & 6 (1074) & 4 (1480) \\
		800 & 73.8 & 27 (499) & 18 (674) & 10 (960) & 8 (1376) & 6 (2025) & 4 (2928)
	\end{tabu}
\end{table}

We now verify that the DtN coarse space provides an approach which is scalable with respect to the number of subdomains $N$. Table \ref{Table:HomoHelmDtNThreshold4o3} details results for a varying number of square subdomains when using a threshold $\eta_{\text{max}} = k^{4/3}$. As well as seeing the iteration counts staying predominantly constant as we increase $k$, they do also as we increase the number of subdomains $N$ (aside from a small number of slightly larger outliers). Note that, while the size of the coarse space increases as we increase $N$, approximately at a rate proportional to $N^{2/3}$ as shown in Figure \ref{Fig:DependencePlots} (in fact, independent of our choice of $\eta_{\text{max}}$), the number of eigenfunctions required per subdomain decreases with $N$. This means the solution of each eigenproblem is much cheaper for large $N$ as they are of smaller size and we require fewer eigenfunctions.

\begin{table}[t]
	\centering
	\caption{Preconditioned GMRES iteration counts when using the two-level method with threshold $\eta_{\text{max}}=k^{4/3}$ for the DtN coarse space and varying the number of subdomains $N$. A uniform decomposition into $\sqrt{N}\times\sqrt{N}$ square subdomains is used.}
	\label{Table:HomoHelmDtNThreshold4o3}
	\tabulinesep=1mm
	\begin{tabu}{cc|ccccccccccccc}
		& & & & & & & & $N$ & & & & & & \\
		$n_{\text{glob}}$ & $k$ & 4 & 9 & 16 & 25 & 36 & 49 & 64 & 81 & 100 & 121 & 144 & 169 & 196 \\
		\hline
		100 & 18.5 & 6 & 6 & 8 & 6 & 6 & 6 & 6 & 6 & 6 & 6 & 6 & 7 & 7 \\
		200 & 29.3 & 6 & 13 & 8 & 6 & 6 & 17 & 7 & 7 & 7 & 7 & 7 & 7 & 7 \\
		400 & 46.5 & 7 & 11 & 7 & 7 & 7 & 7 & 7 & 7 & 10 & 20 & 7 & 7 & 7 \\
		800 & 73.8 & 7 & 9 & 9 & 7 & 7 & 7 & 7 & 7 & 7 & 7 & 8 & 7 & 7
	\end{tabu}
\end{table}

We now turn our attention to the heterogeneous case. Table \ref{Table:HeterogeneousLayersHelm5x5DtNThreshold4o3} (left) gives results for the alternating layers wave guide problem (see Figure \ref{Fig:WaveGuideLayersA}) for varying angular frequency $\omega$, contrast in wave speed $\rho$, and number of subdomains $N$ when using $\eta_{\text{max}} = k_{j}^{4/3}$ in subdomain $\Omega_{j}$. The picture painted is now rather different from the homogeneous case. While for some choices of $N$ iteration counts remain robust to wave number, in general they degrade as $\omega$ increases. The best results are for $N = 4$, $16$, and $64$ (powers of 2) while the poorest are with large $N$. More generally, if the subdomains are close to being aligned with the jumps in $k$ we obtain better results, otherwise robustness is lost. We note, however, that iteration counts are robust to large contrasts $\rho$. We confirm that the disparate trends observed for the alternating layers problem are due to the geometrical aspects of the problem by considering instead the diagonal layers problem (see Figure \ref{Fig:WaveGuideLayersB}). Results for the diagonal layers problem are given in Table \ref{Table:HeterogeneousLayersHelm5x5DtNThreshold4o3} (right) and now show that any robustness to the wave number is, in general, lost for the heterogeneous problem. We note that increasing the threshold to $\eta_{\text{max}} = k_{j}^{3/2}$ does not improve this assessment. Nonetheless, the DtN approach remains robust to increasing the number of subdomains $N$.

\begin{table}[t]
\centering
\caption{Preconditioned GMRES iteration counts for the heterogeneous layers problem when using the two-level method with threshold $\eta_{\text{max}}=k_{j}^{4/3}$ for the DtN coarse space and varying the number of subdomains $N$. A uniform decomposition into $\sqrt{N}\times\sqrt{N}$ square subdomains is used.}
\label{Table:HeterogeneousLayersHelm5x5DtNThreshold4o3}
\tabulinesep=1mm
\begin{tabu}{ccc|ccccccccccc|ccccccccccc}
	& & & \multicolumn{11}{c|}{Alternating layers problem} & \multicolumn{11}{c}{Diagonal layers problem} \\
	& & & & & & & & $N$ & & & & & & & & & & & $N$ & & & & & \\
	$n_{\text{glob}}$ & $\omega$ & $\rho$ & 4 & 9 & 16 & 25 & 36 & 49 & 64 & 81 & 100 & 121 & 144 & 4 & 9 & 16 & 25 & 36 & 49 & 64 & 81 & 100 & 121 & 144 \\
	\hline
	& & 10 & 6 & 6 & 7 & 6 & 6 & 14 & 6 & 16 & 13 & 17 & 13	& 7 & 12 & 15 & 18 & 18 & 18 & 18 & 18 & 17 & 17 & 18 \\
	100 & 18.5 & 100 & 6 & 6 & 7 & 6 & 6 & 14 & 6 & 16 & 13 & 17 & 13 & 11 & 20 & 21 & 20 & 19 & 19 & 18 & 18 & 17 & 16 & 16 \\
	& & 1000 & 6 & 6 & 7 & 6 & 6 & 14 & 6 & 16 & 13 & 17 & 13 & 11 & 20 & 21 & 20 & 19 & 19 & 18 & 18 & 17 & 16 & 16 \\
	\hline
	& & 10 & 8 & 8 & 9 & 8 & 9 & 29 & 9 & 29 & 25 & 32 & 22 & 9 & 17 & 20 & 19 & 19 & 23 & 21 & 28 & 28 & 28 & 28 \\
	200 & 29.3 & 100 & 8 & 7 & 9 & 8 & 8 & 28 & 8 & 28 & 23 & 30 & 20 & 16 & 28 & 30 & 30 & 28 & 29 & 27 & 27 & 25 & 25 & 24 \\
	& & 1000 & 8 & 7 & 9 & 8 & 8 & 28 & 8 & 28 & 23 & 30 & 20 & 16 & 28 & 30 & 29 & 28 & 29 & 27 & 27 & 25 & 25 & 24 \\
	\hline
	& & 10 & 9 & 7 & 7 & 7 & 7 & 25 & 7 & 29 & 25 & 45 & 22 & 10 & 18 & 20 & 25 & 26 & 26 & 25 & 26 & 26 & 32 & 29 \\
	400 & 46.5 & 100 & 9 & 7 & 7 & 7 & 7 & 24 & 7 & 28 & 24 & 44 & 22 & 22 & 39 & 43 & 43 & 40 & 40 & 39 & 37 & 37 & 39 & 35 \\
	& & 1000 & 9 & 7 & 7 & 7 & 7 & 24 & 7 & 28 & 24 & 44 & 22 & 22 & 38 & 43 & 42 & 40 & 40 & 39 & 37 & 37 & 39 & 34 \\
	\hline
	& & 10 & 8 & 10 & 8 & 11 & 9 & 30 & 8 & 34 & 27 & 36 & 33 & 11 & 19 & 24 & 29 & 28 & 30 & 31 & 34 & 37 & 40 & 41 \\
	800 & 73.8 & 100 & 8 & 10 & 8 & 11 & 9 & 38 & 8 & 43 & 33 & 39 & 31 & 32 & 52 & 62 & 61 & 60 & 58 & 58 & 54 & 53 & 53 & 51 \\
	& & 1000 & 8 & 10 & 8 & 11 & 9 & 38 & 8 & 43 & 33 & 39 & 31 & 32 & 55 & 60 & 60 & 59 & 56 & 56 & 54 & 52 & 51 & 49
\end{tabu}
\vspace*{-1mm}
\end{table}

We now show that the sensitivity of the DtN approach is not solely due to the heterogeneity of the media by reconsidering the homogeneous problem but using non-uniform subdomains, which we compute using METIS. Results for this case are given in Table \ref{Table:HomoHelmDtNThreshold4o3METIS} where we see a slow but definite increase in iteration counts as $k$ increases. Again, we see robustness to the number of subdomains but lose robustness to the wave number. Note that this persists even for $\eta_{\text{max}} = k^{3/2}$. Nonetheless, in our DtN approach we still have rather few GMRES iterations required to compute the solution when $k$ is relatively large (in this case up to $k = 117.2$).

\begin{table}[t]
\centering
\caption{Preconditioned GMRES iteration counts (above) and size of the coarse space (below) when using the two-level method with threshold $\eta_{\text{max}}=k^{4/3}$ for the DtN coarse space and varying the number of subdomains $N$. A non-uniform decomposition into $N$ subdomains is used.}
\label{Table:HomoHelmDtNThreshold4o3METIS}
\tabulinesep=1mm
\begin{tabu}{cc|ccccccccccccc}
	& & & & & & & & $N$ & & & & & & \\
	$n_{\text{glob}}$ & $k$ & 4 & 9 & 16 & 25 & 36 & 49 & 64 & 81 & 100 & 121 & 144 & 169 & 196 \\
	\hline
	100 & 18.5 & 7 & 7 & 11 & 7 & 7 & 7 & 7 & 7 & 7 & 7 & 7 & 7 & 7 \\
	200 & 29.3 & 8 & 9 & 10 & 7 & 11 & 9 & 7 & 7 & 7 & 7 & 7 & 7 & 7 \\
	400 & 46.5 & 8 & 10 & 10 & 13 & 14 & 16 & 9 & 14 & 15 & 9 & 8 & 8 & 8 \\
	800 & 73.8 & 8 & 10 & 12 & 12 & 15 & 15 & 13 & 17 & 11 & 14 & 10 & 12 & 17 \\
	1600 & 117.2 & 10 & 12 & 12 & 15 & 16 & 16 & 17 & 17 & 17 & 16 & 16 & 16 & 16 \\
	\hline
	100 & 18.5 & 75 & 158 & 219 & 303 & 397 & 476 & 558 & 644 & 740 & 829 & 923 & 1024 & 1118 \\
	200 & 29.3 & 135 & 282 & 418 & 558 & 677 & 860 & 1003 & 1123 & 1275 & 1435 & 1588 & 1731 & 1867 \\
	400 & 46.5 & 241 & 516 & 751 & 1001 & 1291 & 1569 & 1818 & 2048 & 2294 & 2596 & 2850 & 3145 & 3366 \\
	800 & 73.8 & 481 & 979 & 1446 & 1919 & 2378 & 2844 & 3261 & 3753 & 4291 & 4651 & 5246 & 5720 & 6126 \\
	1600 & 117.2 & 925 & 1857 & 2639 & 3566 & 4408 & 5244 & 6201 & 7008 & 7909 & 8770 & 9563 & 10448 & 11402
\end{tabu}
\end{table}

\section{Conclusions}
\label{sec:5-Conclusions}

In this work we have investigated a two-level domain decomposition approach to solving the heterogeneous Helmholtz equation. Our focus has been on the Dirichlet-to-Neumann coarse space and how the approach depends on the threshold to select which eigenfunctions go into the coarse space. We have seen that the threshold in \cite{ConenEtAl-DtN} can be improved in order to give wave number independent convergence with only moderate added cost due to the larger coarse space. However, this is only true for the homogeneous problem with sufficiently uniform subdomains. In particular, convergence depends on the wave number for a general heterogeneous problem.

In order to obtain fully wave number independent convergence for Helmholtz problems, a stronger coarse space is needed. A recent approach that achieves this, based on a related GenEO-type method, can be found in \cite{BootlandEtAl-H-GenEO}.

\bibliographystyle{spmpsci}
\bibliography{refs}

\begin{thebibliography}{1}
\providecommand{\url}[1]{{#1}}
\providecommand{\urlprefix}{URL }
\expandafter\ifx\csname urlstyle\endcsname\relax
  \providecommand{\doi}[1]{DOI~\discretionary{}{}{}#1}\else
  \providecommand{\doi}{DOI~\discretionary{}{}{}\begingroup
  \urlstyle{rm}\Url}\fi

\bibitem{BabuskaAndSauter-Pollution}
Babu\v{s}ka, I.M., Sauter, S.A.: Is the pollution effect of the {FEM} avoidable
  for the {Helmholtz} equation considering high wave numbers?
\newblock SIAM J. Numer. Anal. \textbf{34}(6), 2392--2423 (1997)

\bibitem{BootlandEtAl-H-GenEO}
Bootland, N., Dolean, V., Jolivet, P.: A {GenEO-type} coarse space for
  heterogeneous {Helmholtz} problems  (2019).
\newblock In preparation

\bibitem{ConenEtAl-DtN}
Conen, L., Dolean, V., Krause, R., Nataf, F.: A coarse space for heterogeneous
  {Helmholtz} problems based on the {Dirichlet-to-Neumann} operator.
\newblock J. Comput. Appl. Math. \textbf{271}, 83--99 (2014)

\bibitem{DoleanEtAl-Book}
Dolean, V., Jolivet, P., Nataf, F.: An Introduction to Domain Decomposition
  Methods: Algorithms, Theory, and Parallel Implementation, vol. 144.
\newblock SIAM (2015)

\bibitem{ErnstAndGander-Helmholtz}
Ernst, O.G., Gander, M.J.: Why it is difficult to solve {Helmholtz} problems
  with classical iterative methods.
\newblock In: I.G. Graham, T.Y. Hou, O.~Lakkis, R.~Scheichl (eds.) Numerical
  Analysis of Multiscale Problems, pp. 325--363. Springer (2012)

\bibitem{NatafEtAl-DtN}
Nataf, F., Xiang, H., Dolean, V., Spillane, N.: A coarse space construction
  based on local {Dirichlet-to-Neumann} maps.
\newblock SIAM J. Sci. Comput. \textbf{33}(4), 1623--1642 (2011)

\bibitem{SaadAndSchultz-GMRES}
Saad, Y., Schultz, M.H.: {GMRES}: A generalized minimal residual algorithm for
  solving nonsymmetric linear systems.
\newblock SIAM J. Sci. Stat. Comput. \textbf{7}(3), 856--869 (1986)

\bibitem{SpillaneEtAl-GenEO}
Spillane, N., Dolean, V., Hauret, P., Nataf, F., Pechstein, C., Scheichl, R.:
  Abstract robust coarse spaces for systems of {PDEs} via generalized
  eigenproblems in the overlaps.
\newblock Numer. Math. \textbf{126}(4), 741--770 (2014)

\end{thebibliography}

\end{document}